\documentclass{tCON2e}    

\usepackage{epstopdf}
\usepackage{subfigure}
\usepackage[longnamesfirst,sort]{natbib}
\bibpunct[, ]{(}{)}{;}{a}{,}{,}

\theoremstyle{plain}
\newtheorem{theorem}{Theorem}

\newtheorem{lemma}{Lemma}

\theoremstyle{definition}

\newtheorem{example}{Example}

\theoremstyle{remark}
\newtheorem{remark}{Remark}

\begin{document}


\def\exp{\mathrm{exp}}  
\def\dtau{\mathrm{d}\tau}
\def\ds{\mathrm{d}s}

\articletype{MANUSCRIPT}

\title{Stabilization and state trajectory tracking problem for nonlinear \break control systems in the presence of disturbances}

\author{
\name{R. Vrabel\textsuperscript{a}$^{\ast}$\thanks{$^\ast$Corresponding author. Email: robert.vrabel@stuba.sk}
}
\affil{\textsuperscript{a}Slovak University of Technology in Bratislava, Faculty of Materials Science and Technology in Trnava, Bottova 25, 
917 01 Trnava, Slovak Republic;
}
\received{August 2016}
}

\maketitle

\begin{abstract}
In this paper we consider the problem of stabilization and tracking of desired state trajectory for a wide range of nonlinear control problems with disturbances. We present the sufficient conditions for the existence of $C^k$ state feedback controllers and the process of their mathematical designing is described. 
\end{abstract}


\section{Introduction and problem formulation}

We will consider the stabilization and state tracking problem for nonlinear control systems with the disturbances in the general form
\begin{equation}\label{problem formulation}
\dot x=f(x,u,w(t)), \quad t\geq 0,
\end{equation}
where $f: \mathbb{R}^n\times \mathbb{R}^m\times \mathbb{R}^p\rightarrow  \mathbb{R}^n,$ $x=(x_1,\dots,x_n)^T$ is the state vector, $\dot x$ is the time derivative of $x$, $ u=(u_1,\dots,u_m)^T$ is the control input variable manipulated by the controller and $w(t)$ represent the bounded, measurable disturbance inputs that cannot be influenced by control.
We will assume that $f$ is a $C^k$  function  ($k\geq2$) in the variables $(x,u)\in\mathbb{R}^{n}\times\mathbb{R}^{m},$ continuous on $\mathbb{R}^{n}\times\mathbb{R}^{m}\times\mathbb{R}^{p}$ and that for every initial state $x(0)$ and input $u$ there exists a unique solution defined on $[0,\infty).$

\section*{Notation and assumptions:}
In this paper, the following notations will be used:
\begin{itemize}
\item $|\, .\, |$ \ -- \ the Euclidean norm on an $n-$dimensional vector space $\mathbb{R}^n$;
\item $||\, .\, ||_{_F}$ \ -- \ the Frobenius matrix norm (or another norm with the properties of submultiplicativity and compatibility with a vector norm);
\item $(\, .\, )^T$ \ -- \ a vector transpose;
\item $\Delta,$ $\tilde\Delta$ \ -- \ the semi-simple Hurwitz matrices of dimension $n$ in the Jordan canonical form (for simplicity, we will restrict ourselves to the matrices with real eigenvalues);
\item $\exp[\ .\ ]$ \ -- \ the exponential function with base Euler's number;
\item $F(e,v,w(t))=:f(e+x_d,v+u_d,w(t))-f(x_d,u_d,w(t))$ \ -- \ the vector field of error dynamics;
\item $A(t)=:\left.J_eF\right\vert_{(0,0)}$\ -- \ the Jacobian matrix of a vector-valued function $F$ with respect to the variable $e$ and evaluated at $(e,v)=(0,0)$   

$(\left.J_eF\right\vert_{(e,v)=(0,0)}=\left.J_xf\right\vert_{(x,u)=(x_d(t),u_d(t))});$
\item $B(t)=:\left.J_vF\right\vert_{(0,0)}$\ -- \ the Jacobian matrix of a vector-valued function $F$ with respect to the variable $v$ and evaluated at $(e,v)=(0,0)$  

$(\left.J_vF\right\vert_{(e,v)=(0,0)}=\left.J_uf\right\vert_{(x,u)=(x_d(t),u_d(t))});$
\item $r(e,v,w(t))=:F(e,v,w(t))-A(t)e-B(t)v$ \ -- \ the Taylor remainder.
\end{itemize}

A large class of control problems consists of stabilization and tracking control for state trajectory
in the presence of the disturbances.  If we are given a desired state
trajectory  $x=x_d(t),$  $t\geq0,$ satisfying (\ref{problem formulation}) for an input $u=u_d(t),$ the goal is to construct a feedback compensator which locally asymptotically stabilizes the system to this trajectory. Using the transformation $e=x-x_d$ and $v=u-u_d,$ the problem of tracking the desired state trajectory may be reformulated in terms of tracking error dynamics which can be written, in general, as a time--varying control system 
\[ 
\dot e=\dot x-\dot x_d=f(x,u,w(t))-f(x_d,u_d,w(t))
\]
\[
=f(e+x_d,v+u_d,w(t))-f(x_d,u_d,w(t))=:F(e,v,w(t)).
\]
The main goal of this paper is to design a state feedback control law $u^*=u_d+v^*(x-x_d,w(t),t)$  [$v^*=v^*(e,w(t),t)$], such that the solution $x(t)$ of (\ref{problem formulation}) [$e(t)$ of $\dot e=F(e,v,w(t))$] asymptotically tracks the desired state trajectory $x_d(t)$ [$e=0$], in the presence of initial state error $x(0)\neq x_d(0)$ [$e(0)\neq0$] and the disturbances $w(t),$ that is,  
\[
x(t)\rightarrow  x_d(t) \ [e(t)\rightarrow 0] \ \mathrm{as} \ t\rightarrow \infty. 
\]
In some cases, for the general control systems (\ref{problem formulation}) and the general time-varying disturbance inputs $w(t),$  it might not be feasible to achieve asymptotic disturbance rejection but only a disturbance attenuation, formulated as a requirement to achieve ultimate boundedness of the tracking error with a prescribed tolerance
\[
|x(t)-x_d(t)|=|e(t)|\leq\varepsilon, \ \mathrm{for\ all}\ t\geq T,
\]
where $\varepsilon$ is a specified (small) positive constant.

Under the assumption that tracking error remains small, we can linearize this system around its equilibrium state $e=0$ (for $v=0$):
\begin{equation}\label{eq:lin}
\frac{de}{dt}\approx A(t)e+B(t)v,\quad \left. A(t)=J_eF\right\vert_{(0,0)},\quad \left.B(t)=J_vF\right\vert_{(0,0)}. 
\end{equation}
If the matrices $A(t),$  $B(t)$ are constant ones and in the absence of time-varying disturbances, for the controllable linearized error system (\ref{eq:lin}) the classical linear control design techniques provide linear feedback control laws
$v^*=-Ke$ which asymptotically stabilize $e=0$ for the closed-loop system. The problem reduces
to calculating a suitable gain matrix $K$ such that $A-BK$ is Hurwitz stable. Moreover, any of
these feedbacks also locally asymptotically stabilize $e=0$ for the original
nonlinear system as follows from the Grobman-Hartman theorem about the local behavior of dynamical systems in the neighborhood of a hyperbolic equilibrium point (\citep{Perko}, p. 127). The location of the closed-loop eigenvalues in the open left half-plane may be
chosen according to the general principle of obtaining fast convergence to zero of the tracking error with a reasonable control effort.  
On the other hand, if system (\ref{eq:lin}) has unstable uncontrollable eigenvalues, then smooth stabilizability is not possible, not even locally. As usual, the critical case is encountered when the linearization has uncontrollable eigenvalues with zero real part.
 

\section{Analysis of the case $m=n$}
The analysis of the problem we divide into two cases:  $m=n$ and $m<n.$ We will not deal with the case $n<m,$  which means that there are more independent actuators than state vector components, that is, the control problem is redundant. The more challenging topics in this area is the design of local stabilizing control laws for the nonlinear systems with more degrees of freedom than control inputs (Section \ref{mlessn}).

The results regarding the case $m=n$ are formulated in the following theorem.
\begin{theorem}[The case $m=n$]\label{main1} 
Let us consider the nonlinear control system (\ref{problem formulation}) and the pair $(x_d(t),u_d(t)).$ Let for all $t\geq0$
\begin{itemize}
\item[$(A1)$] $|x_d(t)|\leq\gamma_{x_d},$ $|u_d(t)|\leq\gamma_{u_d}$ and $|w(t)|\leq\gamma_w$ for some nonnegative constants $\gamma_{x_d},$ $\gamma_{u_d}$ and $\gamma_w;$
\item[$(A2)$] rank of the matrix $B(t)$ is equal to $n;$
\item[$(A3)$] $||B^{-1}(t)||_{_F}\leq\beta_1\exp[\lambda_* t]$  for some constants $\beta_1>0$ and $\lambda_*\geq0,$
\[
\left(||B^{-1}(t)||_{_F}\leq\beta_1\exp[\lambda_* t]\Longleftrightarrow ||B(t)||_{_F}\geq\frac{\sqrt n}{\beta_1}\exp[-\lambda_* t]  \right);
\]
\item[$(A4)$] $||B^{-1}(t)A(t)||_{_F}\leq\beta_2\exp[\lambda_* t],$ for some constant $\beta_2\geq0;$
\item[$(A5)$] there exist the constants $\alpha>0,$ $\beta_3\geq0$ and $\kappa>0$ such that $|r(e,0,w(t))|\leq\beta_3|e|^\alpha,$ for $|e|\leq\kappa,$  where 
\[
r(e,0,w(t))=f(e+x_d,u_d,w(t))-f(x_d,u_d,w(t))-A(t)e.
\]
\end{itemize}
Then there exists a $C^k$ in the variable $x-x_d$ control law $u^*=u_d(t)+v^*(x-x_d,w(t)),$ which makes the desired trajectory $x_d(t)$ of (\ref{problem formulation}) locally asymptotically (exponentially) stable, that is, $e(t)=x(t)-x_d(t)\rightarrow 0$ for $t\rightarrow \infty$ and $x(0)$ satisfying $|x(0)-x_d(0)|<\delta$ for some $\delta>0.$ 

More concretely, for an arbitrary semi-simple Hurwitz matrix $\Delta$ in Jordan canonical form with the eigenvalues $\lambda_i,$ $i=1,\dots, n$  satisfying
\[
2\max\left\{\frac{\lambda_*}{\alpha},\lambda_*\right\}+\max\limits_{i=1,\dots,n}\{\mathrm{Re}\lambda_i\}<0
\]
we have 
\[
x(t)-x_d(t)=\exp[\Delta t](x(0)-x_d(0))\ \mathrm{for} \ t\geq0.
\]
\end{theorem}

To be the system stabilizable by feedback control laws, the assumption $m=n$ may be justified in some cases or significantly simplifies and accelerates the calculations: 

a) Let us consider a driftless $C^1$ nonholonomic control system 
\[
\dot x=\sum\limits_{i=1}^mg_i(x)u_i
\]
and a desired state trajectory  $x_d(t)=x_d$ for $u_d(t))=0$ with constant vector $x_d\in \mathbb{R}^n.$ For tracking error dynamics we obtain differential equation
\begin{equation}\label{error:nonholonomic}
\dot e=\sum\limits_{i=1}^mg_i(e+x_d)v_i=:F(e,v).
\end{equation}
There exists a stabilizing linear control law $v^*=-Ke$ for $\dot e=F(e,v)$  provided the unstable eigenvalues of 
the linearized system are controllable and there exists no stabilizing control law if the linearized system has an unstable eigenvalue which is 
uncontrollable. The Brockett's topological result adapted to error system (\ref{error:nonholonomic}) states that a necessary condition for the existence of a $C^1$ feedback control $v=v^*(e)$ that makes $e=0$ locally
asymptotically stable is that the image of $F(e,v)=\sum\limits_{i=1}^mg_i(e+x_d)v_i$ contains an open neighborhood of $e=0.$
If the vectors $g_i(e+x_d)$ are linearly independent at $e=0,$ then $m=n$ is a necessary and sufficient condition for $C^1$ stabilizability of error system at $e=0$ (\citep{Brockett}, p. 187). For the system (\ref{error:nonholonomic}) we have $A=(0)$ and if the vectors $g_i(e+x_d)$ are linearly independent at $e=0,$ then the assumption $m=n$ implies controllability of a linearization and the error system (\ref{error:nonholonomic}) is locally asymptotically stabilizable with a linear feedback control $v^*=-Ke.$
Notice that the Brockett's theorem does not apply to time-varying feedback laws of the form $v= v(e,t).$ 
More on the topic of stabilization and trajectory tracking of the nonholonomic systems can be found in \citep{GeWangLee,TiaLi,Shi} or in the book \citep{Jarzebowska}.  

b) For the linear time-invariant (LTI) error systems $\dot e=Ae+Bv$ with invertible matrix $B$ we obtain the simple formula for a calculation of state feedback gain matrix $K,$ namely
\[
A-BK=\Delta\Rightarrow K=B^{-1}(A-\Delta),
\]
where $\Delta$ is preassigned Hurwitz matrix, globally asymptotically stabilizing error dynamics at $e=0.$ The equality $K=B^{-1}(A-\Delta)$ is a LTI version of Theorem~\ref{main1}.

In the last decades the stabilization and tracking problems have been intensively studied. Both problems - stabilization of the nonlinear control systems and state trajectory tracking represent the important class of problems in engineering practice.

For example, in the paper \citep{PanWang} a flatness based robust active disturbance rejection control technique scheme with tracking differentiator is proposed for the problem of stabilization and tracking control of the $X-Z$ inverted pendulum, which is widely used in laboratories to implement and validate new ideas emerging in the control engineering.

The paper \citep{Chipofya} presents a solution to stabilization and trajectory tracking of a quadrotor system using a model predictive controller designed using a special type of orthonormal functions. 

A time-varying adaptive controller at the torque level is designed in the paper \citep{Wang} to simultaneously solve the stabilization and the tracking problem of unicycle mobile robots with unknown dynamic parameters. 

In all of these papers, and many others that have appeared in the literature, the specific problems are investigated and the corresponding techniques are developed. In contrast, in this paper we attempted to state the results as generally as possible. 
Therefore, it is possible, that stronger results can be obtained for some special forms of (\ref{problem formulation}) using other methods. One such case is analyzed at the end of this paper in Remark \ref{remark}. Now we proceed to the proof of Theorem \ref{main1}.
 
\subsection{Proof of  Theorem \ref{main1}}
Let $\Delta$ is as yet an arbitrary Hurwitz matrix in the Jordan canonical form and without loss of generality we will assume that for every eigenvalue $\lambda_i,$ $i=1,\dots,n$ its algebraic multiplicity is equal to the geometric multiplicity, that is, $\Delta$ is a semi-simple matrix. Although $m=n,$ we will keep in the proof the original notation, $n$ and $m.$ 

The tracking error dynamics $\dot e=F(e,v,w(t))$ can be expressed using the Taylor series expansion of $F$ in the form
\begin{equation}\label{Taylor_exp}
\dot e=\left[A(t)e+B(t)v+r(e,v,w(t))-\Delta e\right]+\Delta e,
\end{equation}
with an invertible for all $t\geq0$ matrix $B(t).$  Now we show that for every $t_0\geq 0$ fixed, there exist the open neighborhoods $\Omega_e(t_0)$ of $0\in\mathbb{R}^n,$  $\Omega_v(t_0)$ of $0\in\mathbb{R}^m,$ and the $C^k$ function $v^*(e,w(t_0)):$ $\Omega_e(t_0)\rightarrow \Omega_v(t_0)$ such that
\[
A(t_0)e+B(t_0)v^*(e,w(t_0))+r(e,v^*(e,w(t_0)),w(t_0))-\Delta e=0
\]
for all $e\in\Omega_e(t_0).$ Equating the term inside the square brackets in (\ref{Taylor_exp}) to zero we have
\[
v=-B^{-1}(t_0)\left[ A(t_0)e+r(e,v,w(t_0))-\Delta e \right]=:K_e(v).
\]
Fix any $e$ sufficiently near $e=0.$ Then $K_e(v)$ is a function of $v$ only and we may apply the following contraction mapping lemma:
\begin{lemma}(compare with \citep{Hartman}, p. 404)\label{contraction}
Let $\mathcal{B}_a=\left\{ z\in\mathbb{R}^q:\ |z|< a \right\}$ denotes the open ball of radius $a$ centred on the origin in $\mathbb{R}^q.$ If the function $g:\ \mathcal{B}_a\rightarrow \mathbb{R}^q$ obeys
\begin{itemize}
\item[(H1)] there is constant $\Gamma<1$ such that $|g(z_1)-g(z_2)|<\Gamma|z_1-z_2|$ for all $z_1,z_2\in\mathcal{B}_a;$
\item[(H2)] $|g(0)|<(1-\Gamma)a,$
\end{itemize}
then the equation $z=g(z)$ has exactly one solution $z^*,$ and $z^*\in\mathcal{B}_a$. 
\end{lemma}
Now we check that the assumptions of Lemma \ref{contraction} are satisfied. First observe that 
\[
r(e,v,w(t_0))=F(e,v,w(t_0))-A(t_0)e-B(t_0)v
\]
and so, because $A(t_0), B(t_0)$ are the linear transformations, $r\in C^k$ and
\[
\left.J_vr\right\vert_{(0,0)}=0
\]
due to the fact that the reminder $r$ contains only higher-order terms of $e$ and $v.$ By continuity, we may choose $a(t_0)>0$ sufficiently small that
\begin{equation}\label{jacobi}
||\left.J_vr\right\vert_{(e,v)}||_{_F}\leq\frac{\exp[-\lambda_*t_0]}{\beta_1+\eta},\ \eta>0,
\end{equation}
\begin{equation*}
||B^{-1}(t_0)\left.J_vr\right\vert_{(e,v)}||_{_F}\leq ||B^{-1}(t_0)||_{_F}||\left.J_vr\right\vert_{(e,v)}||_{_F}\leq\frac{\beta_1}{\beta_1+\eta}=:\Gamma<1,
\end{equation*}
and
\[
|K_e(v_1)-K_e(v_2)|\leq ||B^{-1}(t_0)\left.J_vr\right\vert_{(e,v)}||_{_F} |v_1-v_2|\leq \Gamma|v_1-v_2|
\]
whenever $|e|,$ $|v|,$ $|v_1|,$ $|v_2|$ are all smaller than $a(t_0).$ It is important to note that $\Gamma$ is a constant independent of $t.$ Also observe that $K_{e=0}(0)=0$ and so we can choose $a'(t_0)\in(0,a(t_0)),$  so that 
\[
|K_e(0)|<(1-\Gamma)a(t_0)
\]
whenever $|e|<a'(t_0).$

We conclude from contraction mapping lemma that, assuming $B(t_0)$ is invertible, there exist $a(t_0),$ $a'(t_0)>0$ such that, for each $e$ obeying $|e|<a'(t_0)$ the system of equations
\[
A(t_0)e+B(t_0)v+r(e,v,w(t_0))-\Delta e=0
\]
has exactly one solution, $v^*(e,w(t_0)),$ satisfying 
\begin{equation*}
|v^*(e,w(t_0))|<a(t_0),
\end{equation*}
$\mathcal{B}_{a'(t_0)}$ $\subset\Omega_e(t_0)$ and $\mathcal{B}_{a(t_0)}\subset\Omega_v(t_0).$ Because $F$ is $C^k$ in the variables $(e,v)$ also the function $v^*(e,w(t_0))$ is $C^k$ in the variable $e.$ Moreover, $v^*$ can be estimated in the following way:
\[
|v^*(e,w(t_0))|=|K_e(v^*(e,w(t_0)))-K_e(0)+K_e(0)|
\]
\[
\leq\mathrm{Lip}(K_e)|v^*(e,w(t_0))|+|K_e(0)|
\]
\[
\leq\Gamma|v^*(e,w(t_0))|+|-B^{-1}(t_0)[A(t_0)e+r(e,0,w(t_0))-\Delta e]|
\]
\[
\leq\Gamma|v^*(e,w(t_0))|+||B^{-1}(t_0)A(t_0)||_{_F}|e|+||B^{-1}(t_0)||_{_F}\beta_3|e|^\alpha
\]
\[
+||B^{-1}(t_0)||_{_F}||\Delta||_{_F}|e|
\]
where $\mathrm{Lip}(K_e)$ is the Lipschitz constant associated to $K_e.$ Hence, for all $t\geq0$ and $e\in\mathcal{B}_{a'(t)},$  $|e|<\kappa$ is
\[
|v^*(e,w(t))|
\]
\begin{equation}\label{estimate2}
\leq\frac{1}{1-\Gamma}\left(||B^{-1}(t)A(t)||_{_F}|e|+||B^{-1}(t)||_{_F}\beta_3|e|^\alpha
+||B^{-1}(t)||_{_F}||\Delta||_{_F}|e|\right).
\end{equation}
For just defined state feedback control law $v=v^*(e,w(t))$  the closed loop dynamics reduces to
\[
\dot e=\Delta e, \ \mathrm{for\ } |e(0)|<\delta=a'(0).
\]
We want to construct a feedback control law $u^*=u_d+v^*$ such that, if $e(t)\rightarrow 0$ then $v^*(e(t),w(t))\rightarrow 0,$ that is, if $x(t)\rightarrow x_d(t)$ then $u^*(t)\rightarrow u_d(t).$ To ensure the vanishing of $v^*(e(t),w(t))$ for $t\rightarrow \infty,$   taking into account that $\Delta$ is a semi-simple matrix in the Jordan canonical form and for such matrix
\[
|e(t)|\leq||\exp[\Delta t]||_{_F}|e(0)|\leq\mu_0\exp\left[\max\limits_{i=1,\dots,n}\{\mathrm{Re}\lambda_i\} t \right]|e(0)|,\ \mu_0>0
\]
and substituting this into (\ref{estimate2}), it is sufficient to choose the matrix $\Delta$ such that 
\begin{equation}\label{ineq}
\max\left\{\frac{\lambda_*}{\alpha},\lambda_*\right\} +\max\limits_{i=1,\dots,n}\{\mathrm{Re}\lambda_i\}<0.
\end{equation}
Analyzing the inequality (\ref{jacobi}) and the entries of the matrix  $\left.J_vr\right\vert_{(e,v)}$ (in the form of Taylor polynomials plus Lagrange remainders) together with the assumption $(A1),$ we can conclude that $a(t)$ and $a'(t)$ decay no faster than $\exp[-\lambda_*t]$ and thus the values $e(t)$ and $v^*(e(t),w(t))$ are defined and fall into the domains $\mathcal{B}_{a'(t)}$ and $\mathcal{B}_{a(t)},$ respectively, for sufficiently small initial values $|e(0)|<\delta$ if a strengthened version of the inequality (\ref{ineq}) holds, namely that
\begin{equation*}
2\max\left\{\frac{\lambda_*}{\alpha},\lambda_*\right\} +\max\limits_{i=1,\dots,n}\{\mathrm{Re}\lambda_i\}<0.
\end{equation*}

The equivalence in the assumption $(A3)$ follows from the inequality
\[
\sqrt n=||B(t)B^{-1}(t)||_{_F}\leq||B(t)||_{_F}||B^{-1}(t)||_{_F}.
\]
Thus Theorem \ref{main1} is proved.
\section{Analysis of the case $m<n$}\label{mlessn}
Analysis of the case of underactuated control system (\ref{problem formulation}) with $m<n$ is much more complicated than the previous one where $m=n.$
The main problem to be overcome is a non-invertibility of the matrix $B(t).$
Let us consider the augmented problem to the original problem (\ref{problem formulation}), namely
\begin{equation}\label{augmented system}
\dot{\tilde x}=f(\tilde x,\tilde u,w(t))+l_{m+1(t)}\tilde u_{m+1}+\dots+l_{n(t)}\tilde u_{n}=:\tilde f(\tilde x,\tilde u,w(t)),
\end{equation}
where $l_{i(t)},$ $i=m+1,\dots,n,$ $t\geq0$ are the column vectors.  We define the pair $(\tilde x_d,\tilde u_d)$ of the augmented system (\ref{augmented system}) associated with $(x_d,u_d)$ of the original system as
\[
(\tilde x_d,\tilde u_d)=(x_d,(u_d^T,\ 0,\dots, 0)^T),
\]
and the error dynamics
\[
\dot{\tilde e}=\tilde f(\tilde e+\tilde x_d,\tilde v+\tilde u_d,w(t))-\tilde f(\tilde x_d,\tilde u_d,w(t))
\] 
\[
=f\left(\tilde e+x_d,\left.\tilde v_i+u_{d_i}\right\vert_{i=1,\dots,m}, w(t)\right)+\sum\limits_{i=m+1}^n l_{i(t)}\tilde v_i
\]
\[
-f(x_d,u_d,w(t))=:\tilde F(\tilde e,\tilde v,w(t))
\]
with $\tilde e(0)=e(0).$ The vectors $l_{i(t)},$ $i=m+1,\dots,n$ are selected such that
\[
\tilde B(t)=\left(B(t)\ \vdots\  l_{m+1(t)}\ \vdots \ \dots\ \vdots \ l_{n(t)}\ \right)
\]
is an invertible matrix. On the basis of the proof of Theorem \ref{main1} there exists a feedback control law $\tilde v^*(\tilde e, w(t))$ which locally asymptotically stabilizes the error dynamics of the augmented system at $(0,0)$ associated with the pair $(\tilde x_d,\tilde u_d).$ Now we can define the feedback control law $v^*=(v^*_1,\dots,v^*_m)^T$ for an error dynamics of the original system as
\begin{equation}\label{association}
v_j^*(e,w(t),t)=\tilde v_j^*(e,w(t)), \ j=1,\dots,m
\end{equation}
where the terms $\tilde v^*_j(\tilde e,w(t)),$ $j=m+1,\dots,n,$ possibly occurring in the argument list of the functions on the right side of (\ref{association}) are replaced by 
\[
\tilde v^*_j(t)=\tilde v^*_j\left(\exp[\tilde\Delta t]e(0)\right),
\]
where the matrix $\tilde\Delta$ satisfies 
\begin{equation}\label{ineq2tilded}
2\max\left\{\frac{\tilde\lambda_*}{\tilde\alpha},\tilde\lambda_*\right\}+\max\limits_{i=1,\dots,n}\{\mathrm{Re}\tilde\lambda_i\}<0.
\end{equation}
Let $E=:\tilde e-e,$ and calculate the difference $\dot E$ between the vector fields for $\dot{\tilde e}(t)$ and  $\dot e(t):$
\[
\dot E=\dot{\tilde e}-\dot e=\tilde F(\tilde e,\tilde v^*(\tilde e),w(t))-F(e,\tilde v^*(e,t),w(t))
\]
\[
=\tilde F(\tilde e,\tilde v^*(\tilde e),w(t))-F(\tilde e-E,\tilde v^*(\tilde e-E,t),w(t))=:H(E,\tilde e, w(t),t),
\]
that is,
\begin{equation}\label{errorE}
\dot E=H(E,\tilde e, w(t),t),\ E(0)=0, \ \tilde e(t)=\exp[\tilde\Delta t]e(0),
\end{equation}
\[
H\in C^k\ \mathrm{in}\ E\ \mathrm{and}\ \tilde e, \ H(0,\tilde e,t)=\sum\limits_{i=m+1}^n l_{i(t)}\tilde v^*_i(\tilde e(t),w(t)). 
\]

The preceding considerations may be summarized in the following theorem. 

\begin{theorem}[The case $m<n$]\label{main2}
Let us consider consider the nonlinear control system (\ref{problem formulation}) and the pair $(x_d(t),u_d(t)).$ Let for all $t\geq0$
\begin{itemize}
\item[$(\tilde{A}1)$] $|x_d(t)|\leq\gamma_{x_d},$ $|u_d(t)|\leq\gamma_{u_d}$ and $|w(t)|\leq\gamma_w$ for some nonnegative constants $\gamma_{x_d},$ $\gamma_{u_d}$ and $\gamma_w;$
\item[$(\tilde{A}2)$] rank of the matrix $B(t)$ is equal to $m;$ 
\item[$(\tilde{A}3)$] $||\tilde B^{-1}(t)||_{_F}\leq\tilde\beta_1\exp[\tilde\lambda_* t],$  where
\[
\tilde B(t)=\left(B(t)\ \vdots\  l_{m+1(t)}\ \vdots \ \dots\ \vdots \ l_{n(t)}\ \right),
\]
with $l_{i(t)}$ chosen so to be $\tilde B(t)$ invertible, $\tilde\beta_1>0,$ $\tilde\lambda_*\geq0,$
\[
\left(||\tilde B^{-1}(t)||_{_F}\leq\tilde \beta_1\exp[\tilde \lambda_* t]\Longleftrightarrow ||\tilde B(t)||_{_F}\geq\frac{\sqrt n}{\tilde \beta_1}\exp[-\tilde \lambda_* t]  \right);
\] 
\item[$(\tilde{A}4)$] $||\tilde B^{-1}(t)A(t)||_{_F}\leq\tilde\beta_2\exp[\tilde\lambda_* t],$ for some constant $\tilde\beta_2\geq0;$ 
\item[$(\tilde{A}5)$] there exist the constants $\tilde\alpha>0,$ $\tilde\beta_3\geq0$ and $\tilde\kappa>0$ such that $|r(\tilde e,0,w(t))|\leq\tilde \beta_3|\tilde e|^{\tilde\alpha},$ for $|\tilde e|\leq\tilde\kappa,$  where 
\[
r(\tilde e,0,w(t))=f(\tilde e+x_d,u_d,w(t))-f(x_d,u_d,w(t))-A(t)\tilde e.
\]
\end{itemize}
Then there exists a $C^k$ in the variable $x-x_d$ control law $u^*=u_d(t)+\tilde v^*(x-x_d,w(t),t),$ which makes the desired trajectory $x_d(t)$ of (\ref{problem formulation}) 
\begin{itemize}
\item[(i)]locally stable (in the sense of Lyapunov), if for every $\varepsilon>0$ there exists $\delta(\varepsilon)>0$ such that $|E(t)|<\varepsilon$ for all $t\geq T(\varepsilon)$ if  $|e(0)|<\delta,$ or
\item[(ii)] locally asymptotically (exponentially) stable, if there exists $\delta_0>0$ such that $E(t)\rightarrow 0$ for $t\rightarrow \infty$ and  for all $|e(0)|<\delta_0.$ Here and above, in (i), the function $E(t)$ is a solution of (\ref{errorE}).
\end{itemize}
\end{theorem}

\begin{example}
As an illustrative example, let us consider the control system
\[ 
{\setlength\arraycolsep{2pt}
\begin{array}{rcl}
       \dot x_1&=&x_1w(t)+x_2+u_1              \\
       \dot x_2&=&x_2+u_1          
       
\end{array}\ \qquad
{\setlength\arraycolsep{1pt}
x_d(t)=\left(\begin{array}{c}
  0  \\
  0  
  
\end{array} \right)},\quad {\setlength\arraycolsep{1pt}
u_d(t)=\left(\begin{array}{c}
  0 
 \end{array} \right)}
}
\]
and its error dynamics along $(x_d,u_d)$
\[ 
{\setlength\arraycolsep{2pt}
\begin{array}{rcl}
       \dot e_1&=&e_1w(t)+e_2+v_1              \\
       \dot e_2&=&e_2+v_1 \ .
\end{array}
}
\]
Augmenting the original system vector field with the vector $l_{2(t)}\tilde u_{2}=(1,\ 0)^T\tilde u_2$ we obtain 
\[ 
{\setlength\arraycolsep{2pt}
\begin{array}{rcl}
       \dot{\tilde e}_1&=&\tilde e_1w(t)+\tilde e_2+\tilde v_1+\tilde v_2              \\
       \dot{\tilde e}_2&=&\tilde e_2+\tilde v_1\           
       \end{array}
}
\]
or
\[
{\setlength\arraycolsep{1pt}
\left(
\begin{array}{c}
  \dot{\tilde e}_1    \\
  \dot{\tilde e}_2
\end{array} 
\right)=\left[\left(
\begin{array}{cc}
  w(t) \ & 1   \\
    0  \ & 1 
\end{array} 
\right)
\left(
\begin{array}{c}
  \tilde e_1    \\
  \tilde e_2  
\end{array} 
\right)+
\left(\begin{array}{cc}
  1 \ & 1   \\
  1\ & 0
\end{array} 
\right)
\left(
\begin{array}{c}
  \tilde v_1    \\
  \tilde v_2  
\end{array} 
\right)\right.
}
\]
\[
{\setlength\arraycolsep{1pt}
-\left.\left(
\begin{array}{cc}
  \tilde\lambda_1 \ & 0  \\
  0\ & \tilde\lambda_2
\end{array} 
\right)
\left(
\begin{array}{c}
  \tilde e_1    \\
  \tilde e_2  
\end{array} 
\right)
\right]+
\left(\begin{array}{cc}
  \tilde\lambda_1 \ & 0  \\
  0\ & \tilde\lambda_2
\end{array} 
\right)
\left(
\begin{array}{c}
  \tilde e_1    \\
  \tilde e_2  
\end{array} 
\right).
}
\]
The expression in the square brackets is equal to zero for
\[
\tilde v^*(\tilde e,w(t))={\setlength\arraycolsep{1pt}
\left(
\begin{array}{c}
  -\tilde e_2 +\tilde\lambda_2 \tilde e_2  \\
  -\tilde e_1w(t)+\tilde\lambda_1 \tilde e_1-\tilde\lambda_2 \tilde e_2
\end{array} 
\right).
}
\]
Hence, by (\ref{association}), the state feedback for the original error system is $v_1^*(e,w)=-e_2 +\tilde\lambda_2 e_2$ and we have the following  error dynamics of original and augmented system: 
\[
{\setlength\arraycolsep{1pt}
\dot{\tilde e}=
 \left(
\begin{array}{c}
  \tilde\lambda_1\tilde e_1  \\
  \tilde\lambda_2\tilde e_2 
\end{array} 
\right),\ \
\dot e=
\left(
\begin{array}{c}
  e_1 w(t) +\tilde\lambda_2 e_2 \\
  \tilde\lambda_2 e_2 
\end{array} 
\right).
}
\]
So 
\[
\dot E={\setlength\arraycolsep{1pt}
\left(
\begin{array}{c}
  \tilde\lambda_1\tilde e_1-e_1w(t) -\tilde\lambda_2 e_2 \\
  \tilde\lambda_2(\tilde e_2 -e_2)
\end{array} 
\right)=\left(
\begin{array}{c}
  \tilde\lambda_1\tilde e_1-(\tilde e_1-E_1)w(t) -\tilde\lambda_2 (\tilde e_2-E_2) \\
  \tilde\lambda_2E_2
\end{array} 
\right),  
}
\]
\[
\ E(0)=0, \ \tilde e_i(t)=\exp[\tilde\lambda_it]e_i(0), \  i=1,2,
\]
where, as follows from (\ref{ineq2tilded}) for $\tilde\lambda_*=0$, the eigenvalues $\tilde\lambda_i,$ $i=1,2$ can be the arbitrary numbers lying in the open left-half of the complex plane (we consider the real numbers only). Because $E_2(t)\equiv0,$  we will focus on the differential equation for the first component of $E:$
\[
\dot E_1-w(t)E_1=h_1(t),\ h_1(t)=:\tilde\lambda_1\tilde e_1-w(t)\tilde e_1-\tilde\lambda_2\tilde e_2\rightarrow 0,\ \mathrm{for}\ t\rightarrow \infty
\]
and its solution for $E_1(0)=0$
\begin{equation}\label{solutionE}
E_1(t)=\exp\left[\int\limits_0^t w(s)\ds \right]\int\limits_0^t h_1(\tau)\exp\left[-\int\limits_0^\tau w(s)\ds \right]\dtau.
\end{equation}
Now let $w(t)$ and $\tilde \lambda_1,$ $\tilde\lambda_2$ are such that
\[
\exp\left[-\int\limits_0^t w(s)\ds \right]\rightarrow \infty \ \mathrm{and}\  \frac{h_1(t)}{w(t)}\rightarrow 0\ \mathrm{for}\ t\rightarrow\infty.
\]
Then, using L'Hospital's Rule to evaluate asymptotics of $E_1(t),$ we have
\[
\lim\limits_{t\rightarrow \infty}E_1(t)=\lim\limits_{t\rightarrow \infty}\frac{\int\limits_0^t h_1(\tau)\exp\left[-\int\limits_0^\tau w(s)\ds \right]\dtau}{\exp\left[-\int\limits_0^t w(s)\ds \right]}=-\lim\limits_{t\rightarrow \infty}\frac{h_1(t)}{w(t)}=0,
\]
for all $e(0).$ Thus, on the basis of Theorem \ref{main2} (ii), the desired trajectory $x_d(t)=0$ is locally asymptotically stable solution (even globally, $\delta_0=\infty$) of the original system with $u=u_1^*(x_1,x_2)=-x_2 +\tilde\lambda_2 x_2.$ 

It is worth noting that, for $w(t)\equiv0,$ we obtain from (\ref{solutionE}) local stability (in the sense of Lyapunov) of $x_d$ only, which is in agreement with Brockett's necessary condition for the existence of a $C^1$ closed loop control law of the form $u^*(x)$ that locally asymptotically stabilizes the nonlinear control system to an equilibrium point.
\end{example}

\begin{remark}\label{remark}
As already mentioned above, Theorem \ref{main1} and  Theorem \ref{main2} capture a broad class of nonlinear control problems, therefore, it can be expected that the stronger results can be obtained for some special classes of (\ref{problem formulation}). For example, let us consider two-input chained system $\dot x=(u_1,\ u_2,\ x_2u_1)^T.$ Applying the technique developed in this paper, for $l_3=(0, 0, 1)^T$ we obtain the state feedback locally asymptotically stabilizing augmented system 
\[
\tilde v^*(\tilde e)=
{\setlength\arraycolsep{1pt}
\left(
\begin{array}{c}
  \tilde\lambda_1\tilde e_1   \\
  \tilde\lambda_2\tilde e_2  \\
  -\tilde e_2\tilde v^*_1+ \tilde\lambda_3\tilde e_3 
\end{array} 
\right)=
\left(
\begin{array}{c}
  \tilde\lambda_1\tilde e_1   \\
  \tilde\lambda_2\tilde e_2  \\
  -\tilde\lambda_1\tilde e_1\tilde e_2+ \tilde\lambda_3\tilde e_3 
\end{array} 
\right),\ \tilde e(t)=\exp[\tilde\Delta t]e(0),
}
\]
with the arbitrary real numbers $\tilde\lambda_1,\tilde\lambda_2, \tilde\lambda_3<0$ ($\tilde\lambda_*=0$) and 
\[
\dot E={\setlength\arraycolsep{1pt}
\left(
\begin{array}{c}
  \tilde\lambda_1\tilde e_1   \\
  \tilde\lambda_2\tilde e_2  \\
  \tilde\lambda_3\tilde e_3 
\end{array} 
\right)-
\left(
\begin{array}{c}
  \tilde\lambda_1 e_1   \\
  \tilde\lambda_2 e_2  \\
  \tilde\lambda_1 e_1  e_2
\end{array} 
\right)=
\left(
\begin{array}{c}
  \tilde\lambda_1 E_1   \\
  \tilde\lambda_2 E_2  \\
  \tilde\lambda_3\tilde e_3-\tilde\lambda_1(\tilde e_1-E_1)(\tilde e_2-E_2)
\end{array} 
\right).
}
\]
The initial state $E(0)=0$ implies $E_1=E_2\equiv0$ and
\[
\dot E_3=\tilde\lambda_3\tilde e_3-\tilde\lambda_1\tilde e_1\tilde e_2, \  \tilde e_i(t)=\exp[\tilde\lambda_1 t]e_i(0),\ i=1,2,3.
\]
Integrating this between $0$ and $t$ we have
\[
E_3(t)=\int\limits_0^t \left(\tilde\lambda_3e_3(0)\exp[\tilde\lambda_3\tau]-\tilde\lambda_1e_1(0)e_2(0)\exp[(\tilde\lambda_1+\tilde\lambda_2)\tau]\right)\dtau
\]
\[
\approx\frac{\tilde\lambda_1e_1(0)e_2(0)}{\tilde\lambda_1+\tilde\lambda_2}-e_3(0). 
\]
As follows from Theorem \ref{main2}\ (i),  the feedback control law $u^*(x)=(\tilde\lambda_1x_1,\tilde\lambda_2x_2)^T$  locally stabilizes (in the sense of Lyapunov) the origin for the closed loop system. We can verify it directly from the explicit solution 
\[
x_1(t)=x_1(0)\exp[\tilde\lambda_1t]
\]
\[
x_2(t)=x_2(0)\exp[\tilde\lambda_2t]
\]
\[
x_3(t)=x_3(0)-\frac{\tilde\lambda_1x_1(0)x_2(0)}{\tilde\lambda_1+\tilde\lambda_2}\left(1-\exp\left[(\tilde\lambda_1+\tilde\lambda_2)t\right]\right),
\]
and so the result is in agreement with Theorem \ref{main2} (i). 

On the other side, the Brockett's necessary condition fails to be satisfied for our system as no point of the form $(0,0,\nu),$ $\nu\neq0$  is in the image of $f,$ and therefore, the system under consideration can not be asymptotically stabilized to the origin by using the control laws of the form $u=u^*(x).$ As is presented in the papers \citep{MurrayWalshSastry} and generally for the chained systems in \citep{TeelMurrayWalsh}, the local asymptotic stabilization to the origin is achieved by the time--varying control law with sinusoids
\[
u_1(x,t)=-x_1+x_3\sin t
\]
\[
u_2(x,t)=-x_2-x_3^2\cos t.
\]
\end{remark}
\section{Conclusions}
The problem of interest is that of regulating the tracking error $e(t)$ around zero by attenuating and possibly rejecting the disturbances $w(t).$
The general framework of control design strategy has been proposed to stabilize (asymptotically or in the sense of Lyapunov) nonlinear control system around the desired state trajectory. We have shown that for $m=n$ the control system with the disturbances can be stabilized at the desired state trajectory by applying the $C^k$ state feedback control law of the form $u^*=u_d(t)+v^*(x-x_d,w(t))$ that the desired state trajectory $x_d(t)$ is a locally asymptotically stable solution of the system $\dot x=f(x,u^*,w(t)).$ 

A weaker result was obtained for the underactuated systems ($m<n$), when the local stabilizability of control system (asymptotic or in the sense of Lyapunov) in the neighborhood of desired state trajectory depends on the properties of associated error system (\ref{errorE}). 

Further efforts could be directed toward weakening the conditions imposed on $||B^{-1}(t)||_{_F}$ and $||B^{-1}(t)A(t)||_{_F}$ by using more general error dynamics system $\dot e=\Delta(e),$ $\Delta(0)=0$ instead of the linear system $\dot e=\Delta e$ with the limiting rate of convergence to $e=0,$ which may be insufficient for some nonlinear control systems.

\section*{Acknowledgments}
This publication is the result of implementation of the project "University
Scientific Park:  Campus MTF STU - CAMBO" (26220220179)
supported by the Research \& Development Operational Program funded by ERDF.


\end{document}